\documentclass{amsart}

\usepackage{amssymb}
\usepackage[all]{xy}
\usepackage{hyperref}

\usepackage{enumitem}   

\setlist[enumerate]{itemsep=.2em,topsep=.2em,leftmargin=1.25em,itemindent=2.0em}

\newtheorem{thm}{Theorem}

\newtheorem{lem}[thm]{Lemma}
\newtheorem{cor}[thm]{Corollary}

\newtheorem{prop}[thm]{Proposition}

\theoremstyle{definition}
\newtheorem{defn}[thm]{Definition}

\newtheorem{say}[thm]{}
\newtheorem{exmp}[thm]{Example}

\newtheorem{ques}[thm]{Question}    

\newtheorem*{ack}{Acknowledgments}      

\newtheorem{defn-thm}[thm]{Definition--Theorem}  
\newtheorem{defn-lem}[thm]{Definition--Lemma}  

\newtheorem{assumption}[thm]{Assumption}
\newtheorem{assumptions}[thm]{Assumptions}  
   
\newtheorem{comment}[thm]{Comment}

\theoremstyle{remark}

\renewcommand{\c}[0]{{\mathbb C}}  

\renewcommand{\o}[0]{{\mathcal O}} 
\newcommand{\z}[0]{{\mathbb Z}}
\newcommand{\n}[0]{{\mathbb N}}

\renewcommand{\a}[0]{{\mathbb A}}

\newcommand{\p}[0]{{\mathbb P}}

\newcommand{\q}[0]{{\mathbb Q}}
\newcommand{\map}[0]{\dasharrow}
\newcommand{\qtq}[1]{\quad\mbox{#1}\quad}

\newcommand{\pico}[0]{\operatorname{\mathbf{Pic}}^{\circ}}

\newcommand{\red}[0]{\operatorname{red}}

\newcommand{\aut}[0]{\operatorname{Aut}}
\newcommand{\auto}[0]{\operatorname{Aut}^{\circ}}

\newcommand{\jac}[0]{\operatorname{Jac}}

\newcommand{\hilb}[0]{\operatorname{Hilb}}

\newcommand{\onto}[0]{\twoheadrightarrow}

\newcommand{\simq}[0]{\sim_{\q}}

\newcommand{\tsum}[0]{\textstyle{\sum}}

\def\into{\DOTSB\lhook\joinrel\to}

\def\loccoh#1.#2.#3.#4.{H^{#1}_{#2}(#3,#4)}

\DeclareMathAlphabet{\mathchanc}{OT1}{pzc}%
                                {m}{it}

\newcommand{\GL}{\mathrm{GL}}

\usepackage[all]{xy}\xyoption{dvips}

\newcommand{\auts}[0]{\operatorname{\mathbf{Aut}}}

\newcommand{\ner}[0]{\operatorname{N\acute{e}r}}
\newcommand{\nert}[0]{\operatorname{N\acute{e}r}^{\rm tot}}
\newcommand{\nero}[0]{\operatorname{N\acute{e}r}^{\circ}}

\newcommand{\aaa}[0]{\operatorname{\mathbb T}}
\newcommand{\ratsec}[0]{\operatorname{RatSec}}

\newcommand{\autsh}[0]{{\mathchanc{Aut}}} 

\begin{document}
\bibliographystyle{amsalpha}


 \title[Abelian fibrations]{Abelian fiber spaces and their\\Tate-Shafarevich twists with  sections}
 \author{J\'anos Koll\'ar}

 \begin{abstract} Starting with an Abelian fiber space, the aim is to construct a Tate-Shafarevich twist that has a rational section.
   
       \end{abstract}

 \maketitle

 Kodaira showed that to any elliptic surface  $S\to C$ without multiple fibers one can associate its relative Jacobian surface  $\jac(S/C)\to C$, which has the same local structure, and also a section.
 Generalizations to  elliptic $n$-folds and to hyperk\"ahler manifolds are studied  in  \cite{MR1260943, MR1242006, MR1272978, MR1401771, MR1929795, MR1918053,  MR4292177, MR4191257, MR4801611, sacca2025, bkv2025, kim2024neronmodelhigherdimensionallagrangian}; see (\ref{prev.work}). 
 This note  considers  Abelian fiber spaces in general.

 \begin{defn}[Abelian torsors and fiber spaces]\label{ab.tor.fibs.defn}
  An {\it Abelian variety} over a field $k$ is a smooth, 
 proper, connected algebraic group $A_k$ over $k$. 
 An {\it Abelian torsor} $X_k$ is a principal homogenous space under an Abelian variety, which can be given by $A_k=\auts^\circ(X_k)\cong \pico\bigl(\pico (X_k)\bigr)$.

  A {\it  fiber space} is a
  proper morphism $\pi:X\to Z$ between normal, irreducible, algebraic spaces, such that
  $\pi_*\o_X=\o_Z$.
   It is an  {\it Abelian fiber space} if 
         general   fibers  are  Abelian torsors.
 \end{defn}

 \begin{defn}[Twisted forms]\label{et.tw.defn}
   A  {\it twisted form} of  a morphism  $\pi:X\to Z$ is obtained by the following procedure.
   Take an open cover $\amalg_i U_i\to Z$ with  $U_{ij}:=U_i\times_Z U_j$, and set  $X_i:=X\times_ZU_i$ and
   $X_{ij}:=X\times_ZU_{ij}$. Choose
   isomorphisms (of $U_{ij}$-schemes)
   $$
   \phi_{ij}:X_{ij}\cong X_{ji},\qtq{such that} \phi_{jk}\circ \phi_{ij}=\phi_{ik}.
   \eqno{(\ref{et.tw.defn}.1)}
   $$
   Use these to patch the $X_i$ together to get a morphism
   $\pi^\Phi:X^\Phi\to Z$.

   Over $\c$ we use the Euclidean topology, in general the  \'etale topology.
   These are also called {\it \'etale} or {\it Tate-Shafarevich}  twists.
   
   A twisted form of an Abelian fiber space is also an
   Abelian fiber space.
        \end{defn}

Our aim is to  prove that, under some natural conditions,
        Abelian fiber spaces have a   twist with a rational section. 
For clarity, we state  first a weaker form over $\c$.
After further definitions, the final version over $\c$ is given in Theorem~\ref{main.thm.v2}. 
The general result  (\ref{main.thm.gen}) over arbitrary base schemes is very similar, but needs more definitions; see (\ref{ab.fibs.defn}) and (\ref{num.sec.defn}).

\begin{assumption}\label{assum.i.a}   Let $\pi:X\to Z$ be an Abelian fiber space   over $\c$, such that 
  $X$ has log terminal singularities,  $K_{X}$ is numerically $\pi$-trivial, $\pi$ is locally projective (in the Euclidean or   \'etale topology) and  of pure relative dimension $d$.
  
   Note that in characteristic 0, any Abelian fiber space is birational to some  $\pi:X\to Z$ such that   $\pi$ is projective, $X$ has  terminal singularities, and $K_{X}$ is numerically $\pi$-trivial; see (\ref{ab.fibs.defn}.3$^\prime$) for details. So these 3 are mild assumptions. (Note that terminal implies log terminal.)
  Having {\em pure} relative dimension holds if $X$ is hyperk\"ahler or if $\dim Z\leq 2$, but not in general.
 \end{assumption}

\begin{thm}[Untwisting of Abelian fiber spaces]\label{main.thm}
  Let $\pi:X\to Z$ be an Abelian fiber space satisfying 
  Assumption~\ref{assum.i.a}. 
  Let   $\eta\in H^{2d}(X, \z)$ be a cohomology class such that $\eta\cap [\mbox{general fiber}]=1$.
         
         Then there is a canonical construction yielding a twisted form
         $$
         \pi_\eta:\aaa_\eta(X/Z)\to Z
         \qtq{with a rational section} s_\eta: Z\map \aaa_\eta(X/Z).
         \eqno{(\ref{main.thm}.1)}
         $$
                 The construction commutes with smooth base changes.
\end{thm}

For fixed $X, Z$, the different $\aaa_\eta(X/Z)$ are birational to each other, but for $\dim Z\geq 2$ one can get infinitely many different ones up to isomorphism, see  (\ref{ell.deg.secs.exmp.1}).

The existence of such $\eta$ is the higher dimensional version of having no multiple fibers. It is an almost necessary  condition, see  Comment~\ref{eta.coment} for a precise statement.

\begin{say}[Plan of the proof]\label{plan.of.pf}
  We  construct $\aaa_\eta(X/Z)$ locally on $Z$, and  patch the  pieces together. We start with the latter.

 (\ref{plan.of.pf}.1) {\it Patching.}   Assume that $Z$ has an open cover $Z=\cup_i U_i$, and over each $U_i$ we have a rational section  $s_i: U_i\map X_i:=\pi^{-1}U_i$. As in (\ref{et.tw.defn}.1),
  the patching is  given by isomorphisms 
  $\phi_{ij}:X_{ij}\cong X_{ji}$ 
  such that $s_j=\phi_{ij}\circ s_i$ on $U_{ij}$.

  If the $\phi_{ij}$  are unique, then the cocyle conditions  $\phi_{jk}\circ \phi_{ij}=\phi_{ik}$ are automatic. How can we achieve this?

  Let $z\in U_i\cap U_j$ be a point such that the fiber $X_z$ is smooth.
  The $\phi_{ij}$ that we aim to construct should restrict to an isomorphism
  $\phi_{ij}|_z: \bigl(X_z, s_i(z)\bigr)\cong  \bigl(X_z, s_j(z)\bigr)$.
  Since $X_z$ is homogeneous,
  such an isomorphism exists, but it is not unique. However, there is a
  unique such  translation-isomorphism, namely translation by $s_j(z)-s_i(z)$.
    Thus the cocyle conditions are  automatic  if we use only translation-isomorphisms for the $\phi_{ij}$.

  (\ref{plan.of.pf}.2) {\it Translation-isomorphisms.}  The notion of translation is more problematic for singular fibers, but in 
  Definitions~\ref{aut.ii.say} and \ref{tors.aut.say} we establish the group of translation-isomorphisms  $\aut^\circ(X/Z)$. 
   Sheafifying in the \'etale topology gives  $\autsh^\circ(X/Z)$.

   (\ref{plan.of.pf}.3) {\it Changing $\eta$.} These consideration also suggest
   that the right place for $\eta$ to live is not  in $H^{2d}(X, \z)$, 
   but in $ H^0(Z, R^{2d}\pi_*\z_X)$. Restriction gives a natural map
   $$
   H^{2d}(X, \z)\to  H^0(Z, R^{2d}\pi_*\z_X).
   $$
   This is an isomorphism if $Z$ is a germ, but otherwise 
   the target can be much larger, so working with
   $ H^0(Z, R^{2d}\pi_*\z_X)$ gives a richer theory.
See also Comment~\ref{eta.coment}.

  The local problem now becomes the following.
  
   (\ref{plan.of.pf}.4)  {\it Local uniqueness.} 
  For a germ $z\in U\subset  Z$ and  $\eta\in H^{2d}(X_U, \z)$, we want  a canonical way to specify a set of
  rational sections  $\ratsec(\eta)=\{s: U\map X_U\}$ that forms a single $\aut^\circ(X_U/U)$-orbit.  We discuss this in more detail in Paragraph~\ref{main.thm.rems}. 

  Once this is done, the local patches automatically glue to give  the global  $\aaa_\eta(X/Z)$.

 (\ref{plan.of.pf}.5) {\it Arbitrary bases.}
We  replace the analytic germs with the strict Henselizations of the local rings, but otherwise the plan is the same.
\end{say}

\begin{defn}[Automorphisms]\label{aut.ii.say}
  For a morphism $\pi:X\to Z$, we let $\aut(X/Z)$ denote the group of automorphism of $X$ that  satisfy $\pi=\pi\circ \phi$.

  By base change, it defines a functor  $V\mapsto \aut(X\times_ZV/V)$. If this functor is represented by a group scheme (or algebraic space), we denote it by
  $\auts(X/Z)$ or  $\auts(X/Z)\to Z$. By a standard Hilbert scheme argument, representability holds if
  $\pi$ is proper and flat \cite[I.1.10]{rc-book}.
  If $\auts(X/Z)\to Z$ is smooth, let $\auts^\circ(X/Z)\subset\auts(X/Z)$ denote the identity component.

  We would like   $\aut^\circ(X/Z)$  to be the global sections of $\auts^\circ(X/Z)$.
   Unfortunately, $\auts(X/Z)$ does not exists in general, and even if it does,
  $\auts(X/Z)$ is frequently not even flat over $Z$.
  
  Nonetheless, if $\pi:X\to Z$ is an Abelian fiber space, in (\ref{tors.aut.say}) we define  $\aut^\circ(X/Z)$ and the corresponding
  sheaf  $\autsh^\circ(X/Z)$.

{\it Note on notation.} The papers
  \cite{abasheva2023shafarevichtategroupsholomorphiclagrangian, kim2024neronmodelhigherdimensionallagrangian}
    define  $\aut^\circ(X/Z)$  for hyperk\"ahler manifolds. The various  definitions most likely agree  when they  all apply.

\end{defn}

\begin{defn}[Translation-automorphisms and twists]\label{trtw.ii.say}
  Let $\pi:X\to Z$ be an Abelian fiber space.
  Since $\aut^\circ(X/Z)$ acts by translations on  general fibers,
we call it the group of
   {\it translation-automorphisms.} It is a commutative group.

  A  twisted form of $\pi$ is     {\it translation-twisted} if
  the patching maps $\phi_{ij}$ in (\ref{et.tw.defn}.1) are
  translation-automorphisms.
  \end{defn}

The stronger  version of  Theorem~\ref{main.thm} is the following. 

\begin{thm}\label{main.thm.v2}
  Let $\pi:X\to Z$ be an Abelian fiber space satisfying 
  Assumption~\ref{assum.i.a}. 
  Let   $\eta\in H^0(Z, R^{2d}\pi_*\z_X)$ be a cohomology class such that $\eta\cap [\mbox{general fiber}]=1$.
         
  Then  $\pi_\eta:\aaa_\eta(X/Z)\to Z$ is a translation-twisted form of $X$, 
        with a rational section $s_\eta: Z\map \aaa_\eta(X/Z)$.
\end{thm}

\begin{comment}\label{eta.coment}
  If $X$ is smooth, then the existence of such an
  $\eta\in H^0(Z, R^{2d}\pi_*\z_X)$ is necessary for $\aaa_\eta(X/Z)$  to exist.

  Indeed, then  $\aaa_\eta(X/Z)$ is also smooth, so the closure of
  $s_\eta(Z)$ defines a cohomology class $\eta'\in H^{2d}\bigl(\aaa_\eta(X/Z), \z\bigr)$ such that $\eta'\cap [\mbox{general fiber}]=1$.

  We see in (\ref{preserve.translate.say}) that $R^{2d}\pi_*\z$ is  the same for $X$ and for
  $\aaa_\eta(X/Z)$. Thus $\eta'$ is identified with  an 
$\eta\in H^0(Z, R^{2d}\pi_*\z_X)$  such that $\eta\cap [\mbox{general fiber}]=1$.

However, if $X$ is singular, then $\aaa_\eta(X/Z)$  sometimes exists even if there is no such $\eta$.
  \end{comment}

\begin{say}[Local case of Theorem~\ref{main.thm}] \label{main.thm.rems} Unfortunately I do not have an a priori definition of what the `canonically constructed' $\aaa_\eta(X/Z)$ should be.
  Instead, we give a procedure to construct it, which goes as follows.

   (\ref{main.thm.rems}.1)  {\it General fibers.}
Using only translation-isomorphisms
  gives the following    over an open, dense set:
  
  Let $Z^\circ\subset Z$ be a dense, open set such that
      $\pi^\circ:X^\circ\to Z^\circ$ is smooth, hence its fibers are
      Abelian torsors.  We take  
      $\aaa_\eta(X/Z)|_{Z^\circ}:=\auts^\circ(X^\circ/Z^\circ)$\footnote{In many papers this is called the
   Albanese fibration of $X^\circ\to Z^\circ$.
   I do not follow this terminology, since
   it is not consistent with \cite[VI.3.3]{FGA}, according to which the Albanese variety is the target of the universal map to an Abelian {\em torsor.}
   Thus, by Grothendieck's definition,  the Albanese  of $X^\circ\to Z^\circ$
   is itself.}.
      This is independent of $\eta$, so the different $\aaa_\eta(X/Z)$ are birational to each other.

      (\ref{main.thm.rems}.2)  {\it Codimension 1 fibers.}  Assume for simplicity that  $X$ is smooth, and $(z,Z)$ is a smooth curve germ.
Then any  section $s:Z\to X$ is regular, and  $s(z)$ is a smooth point of the fiber $X_z$. Thus, in order to specify $\aaa_\eta(X/Z)$, 
we need to choose a local section $s:Z\to X$. Which one?

For minimal elliptic surfaces the local analytic automorphism group acts transitively on the set of local sections, hence all choices lead to isomorphic
objects.
(This is why $\eta$ does not play a role  for elliptic surfaces.)  However, this fails in higher dimensions; see for example \cite[21]{k-neron}.  The best one can say is that by \cite[Thm.1]{k-neron}, the local analytic automorphism group acts transitively on the set of local sections that pass through the same irreducible component of the fiber $X_z$.

Thus, if  $X$ is smooth and $(z,Z)$ is a smooth curve germ, then a choice of a section $s:Z\to X$ up to $\aut^\circ(X/Z)$
 is equivalent to a choice of a reduced
 irreducible component of the central fiber $X_z$.
 If $X$ is singular, it is better to view this as choosing an
irreducible component of the central fiber of the N\'eron model of $X\to Z$ (\ref{ner.mod.say}). 

If  $X$ is arbitrary and $\dim Z\geq 2$,   we need a consistent way to pick an
irreducible component of each fiber of the N\'eron model,   which is defined   outside a codimension $\geq 2$ subset $Z_2\subset Z$.

It turns out that one can use $\eta$ to do this (in a somewhat convoluted way), but  $\eta$ carries much more information than what was needed so far.

      (\ref{main.thm.rems}.3)  {\it Higher codimension fibers for  elliptic 3-folds.}
      
The papers \cite{r-mmc3, MR1260943, MR1242006, MR1272978, MR1401771, MR1929795, MR1918053} describe the situation for elliptic 3-folds  $X\to Z$.
The previous considerations determine  $\aaa(X/Z)$ except over finitely many points $z_i\in Z$, but there are usually infinitely many ways to extend to  twisted forms $\aaa(X/Z)\to Z$. The different models are related to each other by sequences of flops. Choosing a relatively ample line bundle then
specifies an extension uniquely.  The Chern class of a line bundle is in $H^2(X, \z)$, and the same for our $\eta$. They are actually not the same, but  have the same role.  See (\ref{ell.3f.exmp}--\ref{ell.3f.exmp.2}) for examples.

 (\ref{main.thm.rems}.4)  {\it  Higher codimension fibers in general.}
Minimal model theory would suggest that even for higher fiber dimensions, one should choose a relatively ample line bundle to get a unique model.

More generally, choosing an effective divisor $D$ that is ample on the generic fiber, relative MMP for  $(X,\epsilon D)$ results in a model $\pi^c: (X^c, \epsilon D^c)\to Z$ such that
$D^c$ is $\pi^c$-ample. However, $X^c\to Z$ is usually not a twist of $X\to Z$; sometimes it is not even equidimensional.

This is where the role of $\eta\in H^{2d}(X, \z)$  fully emerges.

Assume for simplicity that $\pi$ is flat, $(z, Z)$ is a germ, and let 
$X_z=\sum_{i\in I} m_iF_i$ be the central fiber.
The assumption $\eta\cap [\mbox{general fiber}]=1$ is equivalent to
$\sum_{i\in I} m_ie_i=1$, where $e_i:=\eta\cap [F_i]\in\z$.

A general local complete intersection of codimension $d$ through a smooth point of $F_i$ gives a  multisection $W_i$ of degree $m_i$ of $X\to Z$.
Intersecting with the generic fiber $X_K$ gives a 0-cycle  $[W_i]_K$ of degree $m_i$.
Thus $\sum_{i\in I} e_i[W_i]_K$ is a  0-cycle   of degree $1$ on $X_K$.
Since $X_K$ is a torsor,  we cannot add points. However, pretending that we can, we get a well-defined point as the `sum' of any 0-cycle   of degree $1$; see 
(\ref{mult.sec.say})  for details.
That is, the formal sum $\sum_{i\in I} e_iW_i$ is transformed into  a  rational section of $X\to Z$.

Going through all possible choices of the $W_i$  leads to a  set of rational  sections, denoted by $\ratsec(\eta)$; see   (\ref{loc.rs.say}).

The key---and quite unexpected---result is that  $\ratsec(\eta)$ forms a single
$\auto(X/Z)$ orbit (\ref{main.thm.loc.sh}).

(\ref{main.thm.rems}.5) {\it Construction of automorphisms.}
The set of rational sections  $\ratsec(\eta)$ is  large. Thus it can form a single
$\auto(X/Z)$ orbit only if $\auto(X/Z)$ is large. So we need to construct many automorphisms.

A rational section $s:Z\map X$ is identified with  a $k(Z)$-point of the generic fiber. Thus, given  rational sections $s_1, s_2:Z\map X$,
translation by $s_2-s_1$ gives  a birational map $\phi_{12}:X\map X$ such that
$s_2=\phi_{12}\circ s_1$.

For minimal elliptic surfaces  $\phi_{12}$ is an automorphism, and we are done.
However,  $\phi_{12}$ is usually not an  automorphism for $\dim X\geq 3$, not even if $X$ is smooth.  In (\ref{main.thm.loc.sh}) we prove that  $\phi_{12}$ is an automorphism whenever
$s_1, s_2$ are in the same $\ratsec(\eta)$.  The arguments rely on 
\cite[Thm.1]{k-neron} over codimension 1 points of $Z$, and on 
\cite[37]{k-neron} over higher codimension  points.

\end{say}

If $\dim Z\geq 2$, then $s_\eta$ is rarely a regular section, but
 (\ref{main.thm.rems}.4) shows the following.

\begin{cor} \label{main.thm.cor}
  Using the notation and assumptions of Theorem~\ref{main.thm},
  assume in addition that every geometric fiber has a unique irreducible component of multiplicity 1. Choose $\eta$ to be 1 on that component, and 0 on the others. Then  $s_\eta: Z\to \aaa_\eta(X/Z)$ is a section. \qed
  \end{cor}

\begin{say}[Properties of $\aaa_\eta(X/Z)$]\label{main.thm.rems.7}
Many properties of $X$ are inherited by twists and translation-twists; these are discussed in Section~\ref{twisted.forms.sect}, see especially (\ref{preserve.translate.say}).

By the above construction
$\aaa_\eta(X/Z)$ is an algebraic space, and $\pi_\eta$ is proper and locally projective.

For the most general version of $\eta$ considered in (\ref{num.sec.loc.defn}) and (\ref{num.sec.defn}),
it can happen that $X$ is projective but $\aaa_\eta(X/Z)$ is not; see 
(\ref{ell.deg.secs.exmp.1}).

It is possible that $X$ and   $\aaa_\eta(X/Z)$ are much closer to each other if $\eta\in H^{2d}(X, \z)$.
\end{say}

\begin{ques}\label{homeo.def.ques}  Using the notation and assumptions of  Theorem~\ref{main.thm}, 
  assume in addition that $X$ is projective and
  $\eta\in H^{2d}(X, \z)$.
  \begin{enumerate}
  \item  Is $\aaa_\eta(X/Z)$ projective?
  \item Is $\aaa_\eta(X/Z)\to Z$ homeomorphic to $X\to Z$?
    \item Is $\aaa_\eta(X/Z)\to Z$ deformation equivalent  to $X\to Z$?
    \end{enumerate}
  \end{ques}

\begin{say}[Previous work]\label{prev.work}
   Kodaira's work on   elliptic  surfaces 
\cite{MR132556, MR184257} was extended to  higher dimensional
  elliptic fibrations in the papers
  \cite{r-mmc3, MR1260943, MR1242006, MR1272978, MR1401771, MR1929795, MR1918053}; see also \cite{k-ell, MR4292177, MR4191257, MR4801611} for other results on elliptic varieties.

  The moduli theory of Abelian varieties has been extensively studied,
   see  \cite{MR0352106, MR1083353,   MR1707764, ale-abvar, k-neron} and the references there.

  Higher dimensional  Abelian  fiber spaces arising as  Lagrangian fibrations of  hyperk\"ahler varieties are  investigated in
   \cite{abasheva2023shafarevichtategroupsholomorphiclagrangian, abasheva2024shafarevichtategroupsholomorphiclagrangian,  sacca2025, bkv2025, kim2024neronmodelhigherdimensionallagrangian,  liu-liu-xu}, all of which give  surveys of the earlier literature.

   Roughly speaking, this paper studies for general Abelian fibrations the questions investigated
   in \cite{abasheva2023shafarevichtategroupsholomorphiclagrangian, abasheva2024shafarevichtategroupsholomorphiclagrangian, bkv2025}, which  aim to construct  deformations that have  rational sections. 


   I believe that if $\pi:X\to Z$ is a Lagrangian fibration of a hyperk\"ahler manifold to which \cite[Thm.C]{abasheva2023shafarevichtategroupsholomorphiclagrangian} or \cite[Thm.1.1]{bkv2025} apply, then 
   $\aaa_\eta(X/Z)$ agrees with the degenerate twistor
   deformation constructed there.
   (The choice of $\eta$ is implicit in the proof of \cite[Thm.3.1]{bkv2025}.)
   However, I do not know when $\aaa_\eta(X/Z)$ is homeomorphic or  deformation equivalent  to $X$ in general; see (\ref{homeo.def.ques}) and (\ref{enr.surf.exmp}).

   \end{say}

\begin{ack}  I thank   I.~Levy and   C.~Xu  for  comments,
  Y.-J.~Kim for an especially long list of corrections and discussions, and
  C.~Voisin for suggesting the terminology `untwisting.'
  
  Partial  financial support    was provided  by  the NSF (grant number
DMS-1901855)  and by the Simons Foundation   (grant number SFI-MPS-MOV-00006719-02).
    \end{ack}

\section{Examples}

As we noted in (\ref{main.thm.rems}.6), they key step is to prove Theorem~\ref{main.thm}  when $Z$ is local. This is in turn equivalent to
identifying a set of rational sections $\ratsec(X,\eta)$ forming a 
 single orbit  under translation-automorphisms.

 The proof gives  a procedure to construct $\ratsec(\eta)$, and in many cases it gives a good  description of the resulting rational sections over codimension 1 points of $Z$.  However, it seems difficult to describe the singularities of the sections over codimension $\geq 2$ points.

 Here we discuss a  few  cases for which a reasonably  complete description of $\ratsec(\eta)$ is known.

\begin{exmp}  Assume that $\pi$ is flat, the fiber $X_z$ has a reduced component $F$, $\eta\cap [F]=1$ and $\eta\cap [\mbox{other components}]=0$.

  Then $\ratsec(\eta)$ consist of those sections $s:Z\to X$ for which
  $s(z)\in F$ is a smooth point of $X_z$.
\end{exmp}

\begin{exmp} \label{dimZ=1.exmps}
  Assume that $X$ is smooth and  $\dim Z=1$.
 By \cite[2.1]{k-neron} the
 N\'eron model is isomorphic to the smooth locus of $X\to Z$.

 If the fiber $X_z=\sum_{i\in I} F_i$ is reduced, then 
           the group of connected components of  $\ner(X/Z)_z$ acts simply transitively on $I$.  Let $[i]$ denote the elements of $I$. Since $\sum_{i\in I}e_i=1$,          $\sum_{i\in I}e_i[i]=[i_\eta]$ for some  $i_\eta\in I$.

          Then  $\ratsec(\eta)$ consists of
          those sections $s:Z\to X$ for which
          $s(z)\in F_{i_\eta}$ is a smooth point of $X_z$.

          The cohomology class of these sections is usually  different from $\eta$.

          If  $X_z$ has some nonreduced components, then again
          $\ratsec(\eta)$ consists of  sections through a reduced component $F_{i_\eta}$, but it is  
less clear how $i_\eta$ is determined.
\end{exmp}

Thus, if $\dim Z=1$ then usually there are infinitely many possibilities for $\eta$, but only finitely many for $\ratsec(\eta)$.
In the next example $\dim Z=2$,  and  infinitely many  choices of $\eta$  lead to different  $\ratsec(\eta)$.  This seems typical for elliptic fiber spaces.

\begin{exmp} \label{ell.3f.exmp} Let  $(z,Z)$ be a surface germ, and $\pi:X\to (z,Z)$  an Abelian fiber space of dimension 3, that is, an elliptic 3-fold. Assume for simplicity that 
  $X, Z$ are smooth, $\pi$ is flat,  and  $K_{X}$ is numerically $\pi$-trivial. Let $X_z=\sum m_iF_i$ denote the fiber.
    Then   $\eta\in H^{2}(X, \z)$ is determined by the values  $e_i:=\eta(F_i)$, and  $\eta\cap [\mbox{general fiber}]=\sum m_ie_i$.

  The construction of $\ratsec(\eta)$  goes as follows.

  Pick  smooth points $x_{ij}\in F_i$,  smooth surface germs 
  $x_{ij}\in W_{ij}$ that are transversal to $F_i$, and signs $\epsilon_{ij}$ such that
  $\sum_j\epsilon_{ij}=e_i$. 

  Then $\o_X\bigl(\sum_{ij} \epsilon_{ij} W_{ij}\bigr)$ is a line bundle on $X$.
  Its degree on the general fibers is $\sum_{ij} \epsilon_{ij} m_i=\sum_i m_ie_i=1$.
  Thus it has a section, whose zero locus consists of a rational section  $W\in \ratsec(\eta)$,  plus possibly some divisors that do not dominate $Z$.
  The set of all such $W$ gives $\ratsec(\eta)$.

  To get a clearer description, note that $W$ is usually not nef.
  Let us run the relative MMP for   $(X, cW)$ for any 
  $0<c\ll 1$. If all fibers  other than $X_z$ are irreducible, then the  MMP steps are flops of curves  in $X_z$. In dimension 3, flops preserve smoothness.
  Thus we end with $(\pi^{\rm m}:X^{\rm m}\to Z, W^{\rm m})$ such that  $W^{\rm m}$ is $\pi^{\rm m}$-nef. Its relative degree is 1, so
  $(W^{\rm m}\cdot F^{\rm m}_\ell)=1$ for one irreducible component
  $F^{\rm m}_\ell\subset F^{\rm m}$, and $0$ for the others.
  Moreover, $F^{\rm m}_\ell$ is a component with multiplicity 1, hence
  $\pi^{\rm m}$ is smooth along smooth points of $F^{\rm m}_\ell$.

  The steps of the MMP are independent of the choice of the $W_{ij}$  by \cite[3.53]{km-book}. Thus we conclude that $\ratsec(\eta)$ consists of
  \begin{enumerate}
  \item sections of  $\pi^{\rm m}$ through smooth points of $F^{\rm m}_\ell$, when viewed on $X^{\rm m}$, and
  \item their birational transforms on  $X$, when viewed on $X$.
  \end{enumerate}
  This description is  clear on   $X^{\rm m}$, but less so on $X$.
  Usually members of  $\ratsec(\eta)$ contain the whole
  central fiber, and are possibly  singular along it.

  Using the  detailed description of elliptic 3-folds  given in \cite{MR1260943, MR1929795}, it should be possible to work out $\ratsec(\eta)$ more explicitly in many cases.
  \end{exmp}

\begin{exmp} \label{ell.deg.secs.exmp.1}
  Let $g_0:=x_0x_1x_2$, $g_1, g_2$ be general cubics in the $x_i$.
Let $z:=(1{:}0{:}0)\in Z\subset \p^2_{\mathbf u}$ be a Zariski open subset,
  and set
  $$
  X:=(\tsum u_ig_i=0)\subset \p^2_{\mathbf x}\times Z.
  $$
  Then $X$ is  smooth and the projection $\pi:X\to Z$
  gives an elliptic fiber space.  Then
  $X_z=(x_0x_1x_2=0)$ consists of 3 lines
  $L_i:=(x_i=0)\subset X_z$, and, for suitable $Z$,  all the other fibers are irreducible.

  By the Lefschetz theorem, $H^2(X, \z)$ is generated by the pull-backs of the hyperplane classes. Thus $\eta'\cap[\mbox{general fiber}]$ is divisible by 3 for every $\eta'\in H^2(X, \z)$.

  However, there are classes  $\eta\in H^0(Z,R^2\pi_*\z_X)$
  for which $\eta\cap[\mbox{general fiber}]=1$.
Indeed,   $\eta$ is uniquely determined by the values
$e_i:=\eta\cap L_i$,   and  $\eta\cap[\mbox{general fiber}]=\sum e_i$.

Note also that the classes $[L_i]\in H_2(X, \z)$  are independent of $i$, but the $\eta\cap L_i$ can be different. 

  The different $\aaa_\eta(X/Z)$ are all isomorphic to each other over
  $Z\setminus\{z\}$, but there are infinitely many non-isomorphic
  $\aaa_\eta(X/Z)$.

   $X\to Z$ has relative Picard number 1, hence so does
  $\aaa_\eta(X/Z)\to Z$.
  By the Matsusaka-Mumford criterion 
  only one of the $\aaa_\eta(X/Z)$ is projective.
  (See  \cite{ma-mu}, \cite[11.39]{k-modbook} or \cite[29]{k-neron}.)
\end{exmp}

In the next example $X$ is a mildly singular elliptic 6-fold, and
$\ratsec(\eta)$ consists of sections that are 
not  $\q$-Cartier.
 Thus it is not  possible to associate a cohomology class in $H^2(X, \q)$ to 
$\ratsec(\eta)$. 

\begin{exmp} \label{ell.3f.exmp.2} Let $L_1\neq  L_2$ be lines in $\p^2$.
  Blow up the point $L_1\cap L_2$ twice on $L_2$  to get exceptional curves
  $E_1, E_2$, and then 2  other points on $L_1$ to get $E_3, E_4$.
  The birational transforms of  $E_1, L_1$ become $(-2)$-curves, hence can be contracted. We get $S$, a degree 5 Del~Pezzo surface with an $A_2$ singularity.
The singular point lies on $E'_2$, but not on $L'_2$.
  Note that 
  $2L'_2+3E'_2\in |-K_S|$, where $C'$ denotes the birational transform of a  curve $C$ on $S$. 

  Thus the universal family of divisors in $|-K_S|$ is a  flat, elliptic fiber space  $\pi:X\to \p^5$, one of whose fibers is $2L'_2+3E'_2$.

  For $i=3,4$ we get rational sections  $W_i\subset X$ defined by $|-K_S|\ni C\mapsto C\cap E'_i$.
  Both are regular at the $2L'_2+3E'_2$ fiber. The $W_i$ pass through the singular point.   $W_i$ is not Cartier, but $3W_i$ is. So itersection numbers make sense, giving that
  $(W_i\cdot E'_2)=\frac13$ and $(W_i\cdot L'_2)=0$. Thus
  $\bigl(W_i\cdot (2L'_2+3E'_2)\bigr)=1$, as it should be.

  In local analytic coordinates, we can write a smoothing of $S$ as
  $Y:=(xy+tz-z^3=0)\subset \a^4$. The local structure of the elliptic fiber space is given by the coordinate projection
  $\rho:   Y\to \a^2_{xt}$. The section corresponding to the $W_i$ is
  $W_Y=(y=z=0)$. Note that $W_Y\subset Y$ is not a $\q$-Cartier divisor.
  It does not seem possible to associate a cohomology class in $H^2(Y, \q)$ to $W_Y$. 
  \end{exmp}

\begin{exmp}\label{enr.surf.exmp}
  Let $S$ be an Enriques surface with K3 cover  $\bar S\to S$.
  Let $\tau_S$ be the corresponding involution on $\bar S$.

  Let $E$ be an elliptic curve and $\tau_E$ a translation of order 2.
  Set
  $$X:=(\bar S\times E)/(\tau_S, \tau_E)\qtq{with projection}\pi:X\to S.
  $$
  Since $\tau_E$ acts trivially on the integral cohomology of $E$, we get that $R^2\pi_*\z_X\cong \z_S$. There is a unique $\eta\in H^0(S, R^2\pi_*\z_X)$ such that $\eta\cap[\mbox{general fiber}]=1$.
  Then   $\aaa_\eta(X)\cong  S\times E$, and $\pi_1(S\times E)\cong \z^2+\z/2$. 

  To compute $\pi_1(X)$, write   $E=\c/\langle 1, \tau\rangle$.
  The universal cover of $X$ is  $\bar S\times \c$, and the group of deck transformations is generated by
  $$
  (\bar s, z)\mapsto (\bar s, z+1)\qtq{and}
  (\bar s, z)\mapsto \bigl(\tau_S(\bar s), z+\tau/2).
  $$
  Thus $\pi_1(X)\cong \z^2$,
  hence $X$ and $\aaa_\eta(X)$ are not homeomorphic.

  They are also distinguished by the differentials
  $$
  H^0(S, R^1\pi_*\z_X)\to H^2(S, \z)\qtq{and}
  H^0(S, R^1\pi_*\z_{S\times E})\to H^2(S, \z).
  $$
  The second one is the 0 map, while the first maps onto the 2-torsion in 
  $H^2(S, \z)$.
  \end{exmp}

 \section{N\'eron models}

  The  notion of N\'eron models over regular, 1-dimensional base schemes was developed in 
  \cite{MR179172}; see \cite{MR861977, blr, k-neron} for treatments.
  The theory is not as satisfactory over higher dimensional bases, but
  there is a natural extension if one ignores points of codimension $\geq 2$. 

  \begin{say}[N\'eron models]\label{ner.mod.say}
  Let $Z$ be a normal scheme with function field $K$.
  Let $A_K$ be an Abelian variety over $K$ and $X_K$ an $A_K$-torsor.

  There is a dense, open subset $Z^\circ$ such that
  $A_K$ extends to an Abelian scheme $A^\circ$ over $Z^\circ$, and
  $X_K$ extends to an  $A^\circ$-torsor  $X^\circ$.

  For each codimension 1 generic point $z\in Z\setminus Z^\circ$ one can construct the N\'eron model over the localization  $\o_{z,Z}$, and then extend them to get the following.

  \smallskip
  {\it Claim \ref{ner.mod.say}.1.} There is a large open subset $Z^*\subset Z$ (that is  $Z\setminus Z^*$  has codimension $\geq 2$), such that over $Z^*$ we have   a smooth group scheme $\ner(A_K)\to Z^*$, called the {\it N\'eron model} of $A_K$, and  
  a $\ner(A_K)$-torsor $\ner(X_K)\to Z^*$, called the {\it N\'eron model} of $X_K$.
  The formation of $\ner(X_K)$ commutes with  \'etale base change. \qed
  \smallskip

  If $X\to Z$ is an Abelian fiber space, then the  N\'eron model of its generic fiber is also denoted by $\ner(X^*/Z^*)$, and called the  N\'eron model of $X\to Z$.

  Note that $\ner(X^*/Z^*)\to Z^*$ is  a smooth group scheme,
but $\ner(X^*/Z^*)$ is  not a 
group scheme over $Z$.

  By definition, any $K$-point  $s_K\in X_K(K)$ extends  to a
  {\em rational\ } section  $s_Z:Z^*\map \ner(X_K)$, which is regular at all codimension 1 points.  Since  $\ner(X^*/Z^*)\to Z^*$ is birational to $X\to Z$, we can
  identify rational sections of
   $\ner(X^*/Z^*)\to Z^*$ with rational sections of $X\to Z$.

  If the Mordell-Weil group is finitely generated, after shrinking $Z^*$ we may assume that every $s_Z$ is regular. However, this choice of $Z^*$ is not preserved by \'etale base changes, so not useful for us.
  \end{say}

  \begin{say}[Total  N\'eron model]\label{tot.ner.mod.say}
    Let  $X_i\to Z$ be 
    $A$-torsors under  a flat, commutative group scheme $A\to Z$. The quotient of $X_1\times_ZX_2$ by the $A$-action
    $$
    \tau_a:  (x_1, x_2)\mapsto  (x_1+a, x_2+(-a))
    $$
    is a third $A$-torsor $X_3$, with a  morphism
    $X_1\times_ZX_2\to X_3$.  This makes the disjoint union of all
    $A$-torsors into an algebraic group scheme.

    Starting with an $A$-torsor $X$, applying this operation $(m-1)$-times we get the $A$-torsor $X^{(m)}$. Thus  $X^{(0)}=A$ and $X^{(1)}=X$.
    The disjoint union 
    $\amalg_{m\in\z} X^{(m)}$ is an algebraic group scheme which sits in an exact sequence
    $$
    0\to A \to \amalg_{m\in\z} X^{(m)} \to \z\to 0.
    $$
    Applying this to the   N\'eron model $\ner(X/Z)$ we get 
         the  {\it total  N\'eron model} of $X/Z$  
         $$
         \nert(X^*/Z^*):=\amalg_{m\in\z} \ner^{(m)}(X^*/Z^*).
         $$
         It is a smooth, $\z$-graded,  algebraic group scheme over $Z^*$, whose
degree 0 component is 
$ \ner^{(0)}(X^*/Z^*)=\ner(A^*/Z^*)$, and its degree 1 component is 
$ \ner^{(1)}(X^*/Z^*)=\ner(X^*/Z^*)$.

The {\it identity component} of the total  N\'eron model---denoted by
$\nero(X^*/Z^*)$---is the  largest open subgroup scheme with connected fibers.
Thus $\nero(X^*/Z^*)=\nero(A^*/Z^*)$.
  \end{say}

  \begin{comment} \label{nero.comm} The proof of  Theorem~\ref{main.thm} simplifies if
     $\nero(X^*/Z^*)$  extends to a  smooth group scheme
    $\nero(X/Z)$ with connected fibers.
    If  $\nero(X/Z)$ exists  then   the action of   $\nero(X/Z)$ on $X$ is regular  by \cite[37]{k-neron}.

    Several steps   of the proof can  be viewed as finding some replacement for (the possibly non-existent)  $\nero(X/Z)$.
    \end{comment}

\begin{defn}[Translation-automorphisms]\label{tors.aut.say}
  Let $X_K$ be an $A_K$-torsor. The only automorphisms of $X_K$ as an  $A_K$-torsor are translations by $A_K$. These form the identity component  $\aut^\circ(X_K)$ of  $\aut(X_K)$.  
  
Based on  this, an automorphism  $\phi$ 
of an Abelian fiber space  $\pi:X\to Z$ is a  {\it translation-automorphism} if
the following properties hold.
\begin{enumerate}
  \item For every $z\in Z$, the restriction $\phi_z\in \aut(X_z)$ is a $k(z)$-point of 
    the identity component $\auts^\circ(X_z)$.
  \item  Since $\phi_K$ is in $A_K(K)$, we get a
     section  $s_\phi:Z^*\map \ner^{(0)}(X^*/Z^*)$. We assume that it lies  in the identity component $ \nero(X^*/Z^*)$.
\end{enumerate}
The group of translation-automorphisms is denoted by
$\aut^\circ(X/Z)$. This defines a sheaf in the \'etale topology of $Z$, which we denote by 
$\autsh^\circ(X/Z)$.

Note that, although some  $\auts^\circ(X_z)$ may not be commutative,
$\aut^\circ(X/Z)$ is commutative, since any group identity can be checked on the generic fiber.
\end{defn}

  \begin{say}[Hilbert scheme of local complete intersections]
    \label{hilb.to.ner.say}
    For a  morphism  $g:X\to Z$
    let $\hilb^{\rm ci}_m(X/Z)$ be the relative Hilbert scheme parametrizing
     0-dimensional, local complete intersection subschemes  $P\subset X$ of  degree $m$ such that $g$ is flat along $P$ with Cohen-Macaulay fibers.
     Then $\hilb^{\rm ci}_m(X/Z)$ is smooth over $Z$, and disjoint union provides  rational multiplication maps
     $$
     \hilb^{\rm ci}_{m_1}(X/Z)\times_Z \hilb^{\rm ci}_{m_2}(X/Z) \map
     \hilb^{\rm ci}_{m_1+m_2}(X/Z).
     $$
     This makes $\hilb^{\rm ci}(X/Z):=\amalg_{m\in \n}\hilb^{\rm ci}_m(X/Z)$
     almost into a smooth semigroup scheme over $Z$, where the
      multiplications maps are defined on dense open subsets only.

   If the  generic fiber $X_K$ is an Abelian torsor, 
   then we  get
       a degree preserving semigroup homomorphism
    $$
    \hilb^{\rm ci}(X/Z)\map \nert(X^*/Z^*).
    \eqno{(\ref{hilb.to.ner.say}.1)}
    $$
     \end{say}

      \begin{say}[Multisections]\label{mult.sec.say}
        Let $\pi:X\to Z$ be a  morphism. A {\it multisection} is a reduced subscheme  $W\subset X$ of pure dimension $=\dim Z$,  such that $\pi|_W:W\to Z$ is finite. If $Z$ is integral then $\deg (\pi|_W)$ is the {\it degree} of  $W$.
        $W$ is a section iff its degree is 1. 

         $W$ is a  {\it rational  multisection} if  $\pi|_W$ is generically finite on each irreducible component.

         For general $z\in Z$, the fiber $W_z$ is a  local complete intersection. This gives a rational map  $h_W:Z\map \hilb^{\rm ci}_m(X/Z)$ for $m=\deg (W/Z)$. Combining it with the map in 
    (\ref{hilb.to.ner.say}.1) we get the following.

         \smallskip
        {\it Claim \ref{mult.sec.say}.1.}  Let $\pi:X\to Z$ be an Abelian fiber space, and
         $W\subset X$  a rational  multisection of degree $m$. Then
         $W$ defines a rational section
         $s_W: Z\map \ner^{(m)}(X^*/Z^*)$. \qed
      \medskip

Since $\nert(X^*/Z^*)$ is a group scheme,
a formal $\z$-linear combination
of rational multisections  $W:=\sum_{i\in I}e_iW_i$  gives
$$
s_{W}:Z\map \ner^{(e)}(X^*/Z^*),\qtq{where}  e=\tsum_{i\in I}e_im_i.
\eqno{(\ref{mult.sec.say}.2)}
$$
In particular, if $\sum_{i\in I}e_im_i=1$, then, using that $X$ is birational to $\ner(X^*/Z^*)$,  we get rational sections 
$$
s_{W}:Z\map \ner(X^*/Z^*),\qtq{equivalently}  s_{W}:Z\map X.
\eqno{(\ref{mult.sec.say}.3)}
$$
 \end{say}

            \begin{defn}[Numerical multisections I]\label{num.sec.loc.defn}
         Let $(z, Z)$ be a normal, local scheme and
         $\pi:X\to Z$ a proper morphism  of pure relative dimension $d$.
         Let $\{F_i:i\in I\}$ be the reduced, irreducible components of
         the geometric central fiber.   
         A {\it numerical multisection} is  a $\z$-linear   map on  geometric   $d$-cycles
         $\eta: \sum m_i[F_i]\mapsto   \sum m_ie_i$, provided
         conjugate $F_i$ have the same $e_i$.

         If $\pi$ is flat, then 
a   numerical multisection  $\eta$ is a  {\it numerical section} if
its value on the   $d$-cycle  $[X_z]$ associated to the central fiber
$X_z$ is 1.

If $\pi$ is not flat, then  one needs to use the
{\it Chow-fiber} $X_z^{\rm ch}$  (\ref{chow.fib.say}) instead.
For now the only thing we need to know about it is that
 $X_z^{\rm ch}=\sum m_i[F_i]$ is a $d$-cycle that is algebraically equivalent to the generic fiber.
If  $\pi$ is flat then $X_z^{\rm ch}=[X_z]$, the cycle associated to the subscheme.

 The following is clear.

\smallskip
         {\it Claim  \ref{num.sec.loc.defn}.1.} Let $X_z^{\rm ch}:=\sum_{i\in I} m_iF_i$ be the  Chow-fiber of $\pi:X\to Z$ and $f_i$ the number of
         geometric irreducible components of $F_i$. Then there is a
         numerical section iff $\gcd(m_if_i: i\in I)=1$. \qed
         
\smallskip
      {\it Remark  \ref{num.sec.loc.defn}.2.}     If we are over $\c$, then a  cohomology class $\eta\in H^{2d}(X, \z)$ defines a   numerical multisection  by setting $\eta(F):=\eta\cap [F]$, which is  a numerical section iff
$\eta\cap[\mbox{general fiber}]=1$.
\end{defn}

            \begin{say}[Construction of local  multisections]\label{loc.ms.say}
              Let $(z,Z)$ be a local, Henselian  scheme, and
         $\pi:X\to Z$ a proper,   equidimensional morphism with central Chow-fiber 
              $X_z^{\rm ch}=\sum_{i\in I} m_iF_i$.  Assume that the $F_i$ are geometrically reduced.

              Let  $x_i\in F_i$ be a point such that  $k(x_i)=k(z)$, and
              $\red X_z$ is smooth at $x_i$.
              Let $x_i\in W_i$ be a  local  complete intersection of
              dimension $=\dim Z$ such that  $X_z^{\rm ch}\cap W_i=\{x_i\}$
(as schemes).              (Any general local  complete intersection has this property.)
Since $(z,Z)$ is Henselian, $W_i\subset X$ is a multisection of degree $m_i$.
            \end{say}

\begin{say}[Construction of $\ratsec(\eta)$]\label{loc.rs.say}
         Let $(z,Z)$ be a local, strictly Henselian  scheme, and
         $\pi:X\to Z$ an  equidimensional, Abelian fiber space with  central Chow-fiber 
         $X_z^{\rm ch}=\sum_{i\in I} m_iF_i$.  Assume that the $F_i$ are geometrically reduced, and let $\eta$ be a numerical multisection. Set $e_i:=\eta(F_i)$ and  $e:=\sum_{i\in I} e_im_i$.

         As in (\ref{loc.ms.say}), 
pick points $x_{ij}\in F_i$ and   multisections
$x_{ij}\in W_{ij}\subset X$  of degree $m_i$.
By (\ref{mult.sec.say}.1) they define rational sections
         $$
         s_{W_{ij}}: Z\map  \ner^{(m_i)}(X^*/Z^*)\subset \nert(X^*/Z^*).
         \eqno{(\ref{loc.rs.say}.1)}
         $$
Choose signs  $\epsilon_{ij}\in \{\pm 1\}$ such that
$\sum_j \epsilon_{ij}=e_i$. 
        By (\ref{mult.sec.say}.2)  $\sum_{ij} \epsilon_{ij}  s_{W_{ij}}$ defines a 
         rational section
         $$
         s(\mathbf W, \eta): Z\map \ner^{(e)} (X^*/Z^*).
         \eqno{(\ref{loc.rs.say}.2)}
         $$
         If  $\eta$ is a numerical section, then $e=1$ and $\ner^{(1)} (X^*/Z^*)=\ner (X^*/Z^*)$ is birational to $X$, hence we can view $s(\mathbf W, \eta)$ as   a 
         rational section
         $$
         s(\mathbf W, \eta): Z\map X.
         \eqno{(\ref{loc.rs.say}.3)}
         $$
We call these {\it $\eta$-generic rational sections} of $X\to Z$. 

Varying the points $x_{ij}$ and the  multisections  $W_{ij}$
gives the set of all  $\eta$-generic rational sections, denoted by
$\ratsec(\eta)$ or  $\ratsec(X, \eta)$.
\smallskip

{\it Remark \ref{loc.rs.say}.4.} The natural choice is to take $\epsilon_{ij}$ to be the sign of $e_i$. The flexbility given by the $\epsilon_{ij}$ will be convenient in the proof of (\ref{main.thm.gen}).   
\end{say}

      A  translation-automorphism $\phi$ (\ref{tors.aut.say}) leaves the $F_i$ invariant, so the points $\phi(x_{ij})$ and multisections
$\phi(W_{ij})$ lead to another rational section in
      $\ratsec(\eta)$. This shows the following.

 \begin{lem}\label{tors.aut.say.3}
           Let $(z,Z)$ be a local, strictly Henselian  scheme, and
           $\pi:X\to Z$ an  equidimensional, Abelian fiber space. Then
           $\ratsec(\eta)$ is $\aut^\circ(X/Z)$-invariant. \qed
           \end{lem}

The following key technical theorem shows that  $\ratsec(X, \eta)$ forms a single orbit under translation-automorphisms of $X\to Z$. The
assumptions (\ref{ab.fibs.defn}.1--6) are mostly mild, and can be ignored for now.

      \begin{thm}\label{main.thm.loc.sh}
Let $(z,Z)$ be a local, strictly Henselian  scheme, 
         $\pi:X\to Z$ an  equidimensional, Abelian fiber space
satisfying the assumptions (\ref{ab.fibs.defn}.1--6), and 
          $\eta$  a numerical section (\ref{num.sec.loc.defn}). 
Then
\begin{enumerate}
  \item The rational sections  $\ratsec(\eta)$ constructed in (\ref{loc.rs.say})
    form a single $\aut^\circ(X/Z)$-orbit.
  \item Given $s_1, s_2\in \ratsec(\eta)$, there is a unique
    translation-automorphism $\phi\in \aut^\circ(X/Z)$ such that
    $s_2=\phi\circ s_1$.
    \end{enumerate}
    
      \end{thm}

      The proof will occupy Sections~\ref{rat.homot.sec}--\ref{non.fl.sec}.

      We have a good understanding of birational automorphisms of $X/Z$, these are the $K$-points of $A_K$. To prove regularity over codimension 1 points of $Z$, we use       the following, which  is a direct combination of
 Theorems 1, 37  and 47 in  \cite{k-neron}.

  \begin{thm}\label{neron.1-37-47.thms}
    Let  $\pi:X\to  Z$ be an Abelian fiber space that satisfies the assumptions
    (\ref{ab.fibs.defn}.1--5).
    Then
    \begin{enumerate}
    \item the action of   $\nero(X^*/Z^*)$ on $X^*$ is regular, and
    \item   global sections  of   $\nero(X^*/Z^*)$ induce small modifications of $X$.   \qed
      \end{enumerate}
  \end{thm}

  \section{Rational  homotopies}\label{rat.homot.sec}

   We say that two rational sections $s_i: Z\map X$  are  {\it $\aut^\circ(X/Z)$-equivalent}
    if there is a $\phi\in \aut^\circ(X/Z)$ such that
    $s_2=\phi\circ s_1$. Thus   (\ref{main.thm.loc.sh}) asserts that sections in
    $\ratsec(\eta)$ are $\aut^\circ(X/Z)$-equivalent to each other.

    We prove this in 2 steps. We define a notion of rational  homotopy of sections, and show in (\ref{sm.homot.to aut0.thm}) that  rational  homotopy
    equivalence implies $\aut^\circ(X/Z)$-equivalence.

The construction of rational  homotopy
equivalence is easy if $\pi:X\to Z$ is flat  (\ref{flat.rhomot.prop}).
The equidimensional  case  uses the same idea, but  there are more details to check; see Section~\ref{non.fl.sec}.

  \begin{defn}\label{smooth.homot.defn}
    Let $\pi:X\to Z$ be an Abelian fiber space, and 
    $s_i: Z\map X$  rational sections for $i=1,2$, which we can view as
    $s_i: Z\map\ner(X^*/Z^*)$.
    A {\it  rational  homotopy} between them is a  diagram
    $$
    \begin{array}{ccc}
      W &  \stackrel{\tau}{\map}  &  \ner(X^*/Z^*) \\
      \sigma_i\uparrow\downarrow p{\ } &&  \pi \downarrow\uparrow s_i \\
      Z & = & Z,
      \end{array}
      $$
    where   $p:W\to Z$ is a smooth surjection with geometrically connected fibers,  $\tau:W\map \ner(X^*/Z^*)$ is a rational map, and
      $\sigma_i: Z\to W$ are regular sections such that
      $s_i=\tau\circ\sigma_i$ (as rational maps).

      Note that $\tau$ is a morphism over the generic point of $Z$, so
      $\tau\circ\sigma_i$ are defined (as  rational maps). 
    \end{defn}

  \begin{thm}\label{sm.homot.to aut0.thm}
    Let $\pi:X\to Z$ be an Abelian fiber space
that satisfies the assumptions
    (\ref{ab.fibs.defn}.1--6).
    Let  $s_i: Z\map X$ be rational sections that  are   rational  homotopy  equivalent to each other.  Then they are  $\aut^\circ(X/Z)$-equivalent.
  \end{thm}

  The unexpected feature is that $\tau$  and  the $s_i$ are only rational maps, yet we claim that  they differ by an automorphism.
  Note also that here  $\pi:X\to Z$ is not assumed equidimensional.
\medskip

  Proof.  As we noted in (\ref{plan.of.pf}.1), the condition $s_2=\rho\circ s_1$
  uniquely defines $\rho$ as a rational map,   it is
  translation by  $s_2-s_1$.

  At a codimension 1 point $z\in Z$, $\tau$ is a morphism by the definition of the N\'eron model.  Since $p^{-1}(z)$ is  geometrically connected,
  $s_1(z)$ and $s_2(z)$ are in the same irreducible component of
  $\ner(X/Z)_z$.
Hence $\rho$ induces a   global section of  $\nero(X^*/Z^*)\to Z^*$, showing (\ref{tors.aut.say}.2).

Regularity of $\rho$  is an \'etale-local question;
  we may thus assume that $Z$ is strictly Henselian.

  We aim to apply \cite[36]{k-neron} to  the rational map
  $$
  \phi: X\times_ZW\map X\times_ZW,
  \eqno{(\ref{sm.homot.to aut0.thm}.1)}
  $$
  defined on the smooth fibers of $\pi$ by
  $$
  \phi: (x, w)\mapsto \bigl(x+\tau(w)-s_1(x),w\bigr).
  $$
  Given any section  $\sigma:Z\to W$, the restriction of $\phi$ to
  $\sigma(Z)$ is a rational, translation morphism $X\map X$. Therefore it 
  is a small modification by (\ref{neron.1-37-47.thms}.2).

  Thus, by \cite[36]{k-neron}, there is a small, projective   modification  $\psi:X\map X'$ such that the composite
  $$
  \phi\circ (\psi, 1_W)^{-1}:  X'\times_Z W\map  X\times_Z W\map X\times_ZW
  \qtq{is an isomorphism.}
  $$
  Note that $\phi$ is an isomorphism over  $\sigma_1(Z)$, so
  $\psi:X\map X'$ must be an isomorphism. Thus
  $\phi$ is an isomorphism.

  Its restriction to  $\sigma_2(Z)$ is translation by  $s_2-s_1$, which is $\rho$. Thus 
  $\rho$  is an  isomorphism.

  Finally, let us restrict (\ref{sm.homot.to aut0.thm}.1) to $\{z\}$ to get an isomorphism
  $ \phi_z: X_z\times W_z\to X_z\times W_z$.
  Since $W_z$ is geometrically irreducible, (\ref{tors.aut.say}.1) also holds.
 \qed
  
 \medskip

   The construction of rational  homotopy
equivalences is easier for flat morphisms.

\begin{prop}\label{flat.rhomot.prop}
  Using the notation and assumptions of (\ref{main.thm.loc.sh}), 
let $s, s'\in \ratsec(\eta)$ be
rational sections. Assume also that $\pi$ is flat.
Then  $s, s'$ are rational  homotopy
equivalent.
\end{prop}

The proof is structured to isolate the step (\ref{flat.rhomot.prop}.1) where flatness is used.
Then in Section~\ref{non.fl.sec} we find a replacement for (\ref{flat.rhomot.prop}.1) for non-flat morphisms.
\medskip

Proof. Let $X_z=\sum m_i F_i$ be the central Chow-fiber.
In (\ref{loc.rs.say}) the sections  $s, s'\in \ratsec(\eta)$ are
constructed using some complete intersection multisections
$W_{ij}, W'_{ij}$ of degrees $m_i$ and signs  $\epsilon_{ij}, \epsilon'_{ij}$.

Note that $\sum_j \epsilon_{ij}=e_i=\sum_j \epsilon'_{ij}$, but
$\sum_j |\epsilon_{ij}|$ and $\sum_j |\epsilon'_{ij}|$ may be different.
Assume that $\sum_j |\epsilon_{ij}|=2r+\sum_j |\epsilon'_{ij}|$
We can make them equal by choosing a new $W^*_i$ and adding
$W^*_i$ to the the set $\{W'_{ij}\}$ $r$-times with a positive sign and 
$r$-times with a negative sign. The section $s'$ is unchanged by this.
Permuting the $W'_{ij}$ as needed, we may thus assume from now on that
$\epsilon_{ij}= \epsilon'_{ij}$ for every $i,j$.

 Since $\pi$ is generically flat,  $W_{ij}$ and $W'_{ij}$ define rational sections
 $\sigma_{ij}, \sigma'_{ij}: Z\map \hilb^{\rm ci}_{m_i}(X/Z)$ (\ref{hilb.to.ner.say}).
 Let us make the following
 \smallskip

 {\it Assumption \ref{flat.rhomot.prop}.1.}  There is a commutative diagram
 $$
    \begin{array}{ccc}
      H_{ij} &  \stackrel{u_{ij}}{\map}  &  \hilb^{\rm ci}_{m_i}(X/Z) \\
      \bar\sigma_{ij}\uparrow \bar\sigma'_{ij} &&   \sigma_{ij}\uparrow \sigma'_{ij} \\
      Z & = & Z,
      \end{array}
    $$
    where  $H_{ij}\to Z$ is 
    smooth with geometrically irreducible fibers, and
    $\bar\sigma_{ij}, \bar\sigma'_{ij}:Z\to H_{ij}$ are morphisms.
    \smallskip

If this holds for every $ij$, we can finish the proof as follows.

 The collection of all the $ W_{ij}$  and $W'_{ij}$  gives   sections
 $\sigma, \sigma'$ of  the multiple fiber product $\times_{ij}H_{ij}$.
 Note that $\times_{ij}H_{ij}\to Z$ is also smooth with geometrically irreducible fibers.

 Composing the $u_{ij}$ with the map in (\ref{hilb.to.ner.say}.1) we get  rational maps
 $\tau_{ij}: H_{ij}\map \ner^{(m_i)}(X^*/Z^*)$. Using that $\sum m_ie_i=1$,
 the fiber product of the $\tau_{ij}$, composed with the group structure of
 $\nert(X^*/Z^*)$,  gives
 $$
 \tau: \times_{ij}H_{ij}\map  \ner(X^*/Z^*).
    \eqno{(\ref{flat.rhomot.prop}.2)}
    $$
 This gives a  rational  homotopy
 equivalence between   $s_1, s_2\in \ratsec(\eta)$.
 \smallskip

 It remains the check that   (\ref{flat.rhomot.prop}.1) holds.

 This is easier if $\pi$ is flat, so assume this from now on.
 For each $i$, let $H_i(z)\subset \hilb^{\rm ci}(X/Z)$ be the subset
 parametrizing length $m_i$  complete intersection subschemes of $X_z$ whose support is contained in $F_i$, and also in the smooth locus of $\red X_z$.
 Then $H_i(z)$ is irreducible.

 Using
 \cite[\href{https://stacks.math.columbia.edu/tag/055W}{Tag
     055W}]{stacks-project}
 there is an open $H_i\subset \hilb^{\rm ci}(X/Z)$ containing
  $H_i(z)$ such that $H_i\to Z$ is smooth with geometrically irreducible fibers.

 Any of the $W_{ij}$ used in (\ref{loc.rs.say})  gives a section
 $\sigma_{ij}$ of $H_i\to Z$, so we can set  $H_{ij}:=H_i$, independent of $j$.
 \qed

\section{Non-flat case}\label{non.fl.sec}

\begin{say}[Chow fiber]\label{chow.fib.say}
  Let    $\pi:X\to Z$ be a proper morphism  of pure relative dimension $d$.
  The {\it Chow-fiber}  $X_z^{\rm ch}$ is a $d$-cycle suported on $X_z$ such that
  \begin{enumerate}
  \item $X_z^{\rm ch}=[X_z]$ if  $\pi$ is flat at the generic points of $X_z$, and
   \item the  $X_z^{\rm ch}$ are algebraically equivalent to each other.
 \end{enumerate}
  It turns out that Chow-fibers exist if $Z$ is normal and
  $\red(X_z)$ is geometrically reduced.
  Here are 2 ways to construct  them.

  \smallskip
      {\it Specialization  \ref{chow.fib.say}.3.}  (See \cite[I.3]{rc-book} for details.)
      
         Let $R$ be the spectrum of a DVR and $g:R\to Z$ a local morphism that maps the generic point of $R$ to the generic point of $Z$.
         Since $\pi$ is not assumed flat, the fiber product $X\times_ZR$ may have torsion along the central fiber. Taking the quotient by this torsion, we get
         $g^{[*]}X\to R$, which is now flat. This assigns multiplicities to the 
         irreducible components of the  geometric central fiber of $g^{[*]}X\to R$.  This defines  $X_z^{\rm ch}$ over a (usually infinite) extension of $k(z)$.
         As long as  $Z$ is normal, the resulting
         multiplicities do not depend on the choice of $g:R\to Z$.
        
         If $\red X_z=\sum_{i\in I} F_i$ is geometrically reduced  (for example, if $k(z)$ is perfect),  then the Chow-fiber exists as a cycle on
         $X$. We assume this from now on and write 
         $X_z^{\rm ch}:=\sum_{i\in I} m_iF_i$.

 \smallskip
     {\it Projection  \ref{chow.fib.say}.4.}   Assume that $(z,Z)$ is local,
     $\pi$ is projective and $k(z)$ is infinite. Choose an embedding  $X\subset \p^N_Z$.
     Let $p:\p^N_Z\map \p^{d+1}_Z$ be a projection whose restriction to $X$ is a
     morphism  $p_X:X\to X'\subset \p^{d+1}_Z$.
     Since $Z$ is normal, $X'$ is a Cartier divisor in $\p^{d+1}_Z$ by the
     Ramanujam-Samuel theorem; see \cite[IV.21.14.1]{ega} or \cite[10.65]{k-modbook} for proofs.  Thus $X'$ is flat over $Z$  with Cohen-Macaulay fibers.

     Let $[X'_z]$ be the cycle associated to the fiber $X'_z$.
     If $\red X_z$ is geometrically reduced, we can choose $p$ such that 
     its restriction  to $\red X_z$ is birational. Then
      $[X'_z]$ can be  identified with a cycle supported on $X_z$, giving the 
     Chow-fiber.

\smallskip
     {\it Remark  \ref{chow.fib.say}.45.}
     Achieving flatness by projecting to a relative hypersurface is the key idea in the notion of {\it K-flatness;} see \cite[Sec.7.1]{k-modbook}.

  \end{say}

In order to prove (\ref{main.thm.loc.sh}), it remains to establish
(\ref{flat.rhomot.prop}) without the flatness assumption.

\begin{prop}\label{nonflat.rhomot.prop}
    Using the notation and assumptions of (\ref{main.thm.loc.sh}), 
let $s, s'\in \ratsec(\eta)$ be
rational sections.   
  Then there is a chain of rational  homotopy
  equivalences connecting $s_1$ and $s_2$.
  \end{prop}

Proof.
As we noted during the proof of (\ref{flat.rhomot.prop}),
we need to check  Assumption~\ref{flat.rhomot.prop}.1.
Instead we show something slightly weaker, which is however enough for the proof to go through.

As before,  the sections  $s, s'\in \ratsec(\eta)$ are
constructed using  multisections
$W_{ij}$ and $W'_{ij}$ of degrees $m_i$.
We work with each pair $W_{ij}, W'_{ij}$ separately, so we drop the subscript.

\medskip
{\it Claim \ref{nonflat.rhomot.prop}.1.} Given $W$ and $W'$, there is a
sequence
$
W=: W^0, \dots,  W^r:=W',
$
such that  (\ref{flat.rhomot.prop}.1) holds for each pair
$ W^c, W^{c+1}$ for $c=0,\dots r-1$.
\medskip

To see this,
let $p_X:X\to X'\subset \p^{d+1}_Z$ be a projection as in (\ref{chow.fib.say}.4). Since $X'\to Z$ is flat, as in (\ref{flat.rhomot.prop}), we can use
$\hilb^{\rm ci}(X'/Z)$ to construct rational  homotopy
equivalences for local complete intersections on $X'$. These then give
 rational  homotopy
equivalences for the pulled-back local complete intersections on $X$.

We are not yet done since usually there is no projection $p_X:X\to X'\subset\p^{d+1}_Z$ such that
both $W$ and $W'$ are pull-backs of  local complete intersections on the same $X'$.
However, we claim that using several projections, we get a
chain of rational  homotopy
equivalences.

Assume that $W$ is an irreducible component of a
complete intersection  $A_1\cap\cdots\cap A_d$. If $B_i$ are general, then
$W$ is also an irreducible component of the
complete intersection  $(A_1+B_1)\cap\cdots\cap (A_d+B_d)$.
Choosing the $B_i$ suitably, we may thus assume that
$W $ (resp.\ $W'$) is an  irreducible component of a
complete intersection   $\cap_i(\ell_i=0)$
(resp.\  $\cap_i(\ell'_i=0)$) where
the $\ell_i, \ell'_i$ are  sections of the same  $\pi$-very ample line bundle $L$.

Fix now the embedding  $X\into \p^N_Z$ given by $L$. From now on we use
fiberwise linear projections $\p^N_Z\map \p^{d+1}_Z$.

First choose a projection whose restriction to $\red X_z$ is an embedding
at both points $F\cap (W\cup W')$. This gives a rational  homotopy
equivalence between some $W_1$ and $W'_1$ such that $F\cap W=F\cap W_1$ and
$F\cap W'=F\cap W'_1$.

We are now reduced to the case when    $F\cap W=F\cap W'$ is the same point $x\in F$.

We have $W=\cap_i(\ell_i=0)$. Let $\ell^*_1$ be a general section of  $L$ vanishing at $x$. Then   $\langle \ell^*_1, \ell_1, \dots, \ell_d\rangle$ defines a projection $\p^N_Z\map \p^{d+1}_Z$  which gives a rational  homotopy
equivalence between 
$$
W=(\ell_1=0)\cap\cdots\cap (\ell_d=0)\qtq{and}
W^1:=(\ell^*_1=0)\cap (\ell_2=0)\cap\cdots\cap (\ell_d=0).
$$
We can iterate the procedure to get a chain of rational  homotopy
equivalences between 
$$
W=(\ell_1=0)\cap\cdots\cap (\ell_d=0)\qtq{and}
W^d=(\ell^*_1=0)\cap\cdots\cap (\ell^*_d=0).
$$
The same argument gives a chain of rational  homotopy
equivalences between $W'$ and $W^d$. \qed

\section{Main theorem}

Trying to get the minimal set  of assumptions for the general version of Theorem~\ref{main.thm.v2} results in a  list that is a bit long. However, the  only strong restriction is
(\ref{ab.fibs.defn}.7). 

       \begin{assumptions}\label{ab.fibs.defn}
  We impose  the following assumptions on
   Abelian fiber spaces   $\pi:X\to  Z$.
   \begin{enumerate}
     \item    $X$  and $Z$ are excellent, normal, integral algebraic spaces, and $K_{X}$ exists.
     \item    $\pi:X\to  Z$ is  locally projective.
       \item  $K_{X}$  is $\q$-Cartier and numerically $\pi$-trivial.
       \item There is a closed subset $Z_2\subset Z$ of codimension $\geq 2$ such that  $X\setminus \pi^{-1}(Z_2)$ has
 log terminal singularities if $Z$ is over a field of characteristic 0, and
 terminal singularities otherwise.
\item There are no $\pi$-exceptional divisors; that is, divisors $D\subset X$ such that $\pi(D)$ has codimension $\geq 2$ in $Z$.
  \item For every $z\in Z$, $\red X_z$ is geometrically reduced.
   \end{enumerate}
As we discuss below, these assumptions are quite mild in
characteristic 0, more generally, when  MMP works.
There is, however a stonger version of
(\ref{ab.fibs.defn}.5) that is needed for the main theorems.
   \begin{enumerate}\setcounter{enumi}{6}
   \item  $\pi$ has pure relative dimension.
      \end{enumerate}
 Note that all the assumptions are preserved by smooth base changes  $U\to Z$.
   
 {\it Comments on the assumptions.}
 
 (\ref{ab.fibs.defn}.1$^\prime$) $K_{X}$ exists if $Z$ is of finite type over a field. More generally, $K_{X}$ exists if  $Z$ is of finite type over an excellent scheme that has  a dualizing complex. See \cite[11.2]{k-modbook} for an overview and
 \cite{lyu-mur} for details.
   
(\ref{ab.fibs.defn}.2$^\prime$) Most likely everything holds if $\pi$ is only proper and also for complex analytic spaces. Local projectivity is used mainly when applying \cite[30 and 36]{k-neron},  and it was also used in the non-flat case (\ref{nonflat.rhomot.prop}). 
   
 (\ref{ab.fibs.defn}.3$^\prime$)
As we already noted after  Assumption~\ref{assum.i.a},
 minimal model theory predicts that any regular
 Abelian fiber space  $Y\to Z$ has a relative minimal model  $\pi^{\rm m}:Y^{\rm m}\to Z$
 with terminal singularities,  such that  $K_{Y^{\rm m}}$  is $\q$-Cartier and
 $\pi^{\rm m}$-nef. The abundance conjecture gives a factorization
 $\pi':Y^{\rm m}\to Z'$ such that $K_{Y^{\rm m}}$  is
 numerically $\pi'$-trivial and $Z'\to Z$ is birational.
From this point of view, (\ref{ab.fibs.defn}.3)  is not restrictive.
 
 Roughly speaking, minimal models are known to exist
        in characteristic 0,  and are conjectured to exist in general;
       the main references are
       \cite{MR2359343, bchm, lyu-mur},  see also \cite{fujino-ssmmp, lai-2009} for Abelian fiber spaces and \cite{sacca2025} for the hyperk\"ahler case.

       Note also that by \cite[9]{k-ell}, if $Z$ is $\q$-factorial and
       $K_{X}$  is $\pi$-nef and numerically $\pi$-trivial over a large open subset $Z^*\subset Z$, then it is numerically $\pi$-trivial everywhere.

  (\ref{ab.fibs.defn}.4$^\prime$)  Most likely in (\ref{ab.fibs.defn}.4) one can allow log terminal singularities in all cases. The terminal assumption is dictated by 
   \cite[Thm.1]{k-neron}. It is  extended to log terminal singularities in \cite[47]{k-neron} using MMP, but currently MMP is known only in  characteristic 0.

   Note also that restrictions on the singularities of $X$ impose similar restrictions on the singularities of $Z$. For example, if $X$ has    log terminal singularities, then
   $Z$ also has  log terminal singularities outside a subset of codimension $\geq 3$, and   potentially log terminal singularities everywhere, at least in  characteristic 0; see  \cite{naka88, MR2153078}.

    (\ref{ab.fibs.defn}.5$^\prime$)  Once (\ref{ab.fibs.defn}.3) holds and $X$ has log terminal singularities,
   one can contract all $\pi$-exceptional divisors to achieve (\ref{ab.fibs.defn}.5); see \cite[8]{k-ell}.

   (\ref{ab.fibs.defn}.6$^\prime$) Note that  (\ref{ab.fibs.defn}.6)  always holds in characteristic 0. It is possible that it is sufficient to require (\ref{ab.fibs.defn}.6) only for closed points $z\in Z$.

   (\ref{ab.fibs.defn}.7$^\prime$)
    Note that  (\ref{ab.fibs.defn}.5)  and (\ref{ab.fibs.defn}.7) are equivalent  if $\dim Z\leq 2$, but having  pure relative dimension is more restrictive in general.
  
  It can  be achieved by suitably blowing up $Z$, though it is not clear that one can get both log terminal singularities and  pure relative dimension           this way.

  The conclusion of the main theorems  might hold
  even if $\pi$ is not pure dimensional.
       \end{assumptions}

  We also need to extend (\ref{num.sec.loc.defn}) to find a replacement for $ H^0(Z, R^{2d}\pi_*\z_X)$ in   (\ref{plan.of.pf}.3). 
     
          \begin{defn}[Numerical multisection II]\label{num.sec.defn}
         Let $\pi:X\to Z$ be a proper morphism of relative dimension $d$.   For a geometric point $z\to Z$, a $d$-cyle on $X_z$ is called a
         {\it vertical  $d$-cycle} on $X$. 

        A $\z$-linear function
        $\eta$ from vertical  $d$-cycles to $\z$  is a
        {\it numerical multisection} if the following hold.
        \begin{enumerate}
          \item   Let $C$ be a regular, irreducible, 1-dimensional scheme, $C\to Z$ a morphism and
        $S\subset X\times_ZC$ an irreducible subscheme of dimension $d+1$ with projection $\pi_S:S\to C$. Then   $c\mapsto \eta\bigl[\pi_S^{-1}(c)\bigr]$ is constant.
         \item  $\eta(F_1)=\eta(F_2)$ if $F_1, F_2$ are Galois conjugate cycles.
        \end{enumerate} 
        A   numerical multisection  $\eta$ is a  {\it numerical section} if          $\eta(\mbox{general fiber})=1$.

        The main examples are the following.
         \begin{enumerate}\setcounter{enumi}{2}
         \item  If we are over $\c$, then any
$\eta\in  H^0(Z, R^{2d}\pi_*\z_X)$
           defines a   numerical multisection  by setting $\eta(F):=\eta\cap [F]$.
        \item If $X$ is regular then any  multisection $W\subset X$ defines a   numerical multisection  by setting $\eta_W(F):=(W\cdot F)$.
          It is a numerical section iff $W$ is a rational section.
        \item Let $X'\to Z$ be a translation twist of $X\to Z$, and $\eta$ a 
numerical (multi)section on $X\to Z$. This gives a $\z$-linear function
$\eta'$ on vertical  $d$-cycles on $X'$ which satisfies (\ref{num.sec.defn}.1),
but need not satisfy (\ref{num.sec.defn}.2).
                              \item  The assumption (\ref{num.sec.defn}.1) considers only  algebraic equivalences where all intermediate cycles are also vertical. 
           Thus it can happen that $Z_1, Z_2$ are numerically equivalent, but
           $\eta(Z_1)\neq \eta(Z_2)$; see (\ref{ell.deg.secs.exmp.1}).
           \item We can also think of $\eta$ as a function on the irreducible components $F$ of fibers over closed points of $Z$, provided $\eta(F)$ is divisible by the number of geometric irreducible components of $F$. 
 \end{enumerate}
          \end{defn}

We can now state the general version of Theorems~\ref{main.thm} and \ref{main.thm.v2}.

          \begin{thm}\label{main.thm.gen} Let $\pi:X\to Z$ be an Abelian fiber space  satisfying the assumptions (\ref{ab.fibs.defn}.1--7), and
            let $\eta$ be a numerical multisection (\ref{num.sec.defn}). 
         Then there is a canonical construction yielding a translation-twisted form
         $$
         \pi_\eta:\aaa_\eta(X/Z)\to Z
         \qtq{with a rational section} s_\eta: Z\map \aaa_\eta(X/Z).
         \eqno{(\ref{main.thm.gen}.1)}
         $$
     The construction commutes with smooth base change and  localization.
\end{thm}

          Proof. For a point $z\in Z$ with  strict Henselization
          $Z^{\rm sh}_z$,
          (\ref{num.sec.defn}.2) guarantees that $\eta$ gives  a well defined
          numerical section over $Z^{\rm sh}_z$. Then 
          (\ref{main.thm.loc.sh}) constructs
          $$
          \aaa_\eta\bigl(X\times_ZZ^{\rm sh}_z/Z^{\rm sh}_z\bigr)
          \qtq{with a rational section.}
          \eqno{(\ref{main.thm.gen}.2)}
          $$
          Since $X\to Z$ is of finite type, these are defined over an \'etale neighborhood of $z\in Z$. Thus we get an \'etale cover
          $U:=\amalg_i (u_i, U_i)\onto Z$ such that we have
          $$
          \aaa_\eta\bigl(X\times_ZU/U\bigr)
          \qtq{and} s_U:U\map \aaa_\eta\bigl(X\times_ZU/U\bigr).
          \eqno{(\ref{main.thm.gen}.3)}
          $$
          So far we know that this gives the right answer over the
          strict Henselization of each $(u_i\in U_i)$. We need to check that
          the same holds over all points.

          Once this is established, the uniqueness of the translation-isomorphisms (\ref{main.thm.loc.sh}.2)  implies that (\ref{main.thm.gen}.3) descends to give (\ref{main.thm.gen}.1).

          Dropping the subscript $i$,  assume that we have
          $\pi_U: X_U\to U$ with fiber $\red X_u=\cup_j F_j$.
          By shrinking $U$ if necessary, we may assume that every other fiber $ X_{u'}$ specializes to $X_u$.  That is, there is an irreducible curve $C\to U$
          with points $ c, c'\in C $ mapping to $ u, u'\in U$
         such that, for every irreducible component
          $Y_m\subset X\times_UC$, the fiber $F'_m$ of $Y_m\to Z$ over $c'$ is irreducible and generically reduced. If this holds then  $\red X_{c'}=\cup_m F'_m$.  Let $\sum_j c_{jm}F_j$ be the specialization of $F'_m$.
          
          The construction of the rational section in (\ref{loc.rs.say}) uses
          local complete intersections $W_j$ meeting $F_j$ transversally. Here we have a free choice, so we may assume that the  $W_j$ also meet
          the $F'_m$  transversally.

          Then we consider  $\sum_j \eta(F_j)W_j$. 
          Each $W_j$ meets $F'_m$ at $c_{jm}$ points, so on $F'_m$, we have
          $\sum_j c_{jm}\eta(F_j)$ transversal intersection points.
          This sum equals  $\eta(F'_m)$.  Thus  $\sum_j \eta(F_j)W_j$ also computes a section in $\ratsec(\eta)$ over $u'$. \qed

          \begin{comment}\label{eta.coment.gen}
            As in  Comment~\ref{eta.coment}, if $X$ is a smooth variety over an algebraicaly closed field, then the existence of a numerical section is 
            necessary for $\aaa_\eta(X/Z)$  to exist; see (\ref{num.sec.defn}.5).
            \end{comment}

    \section{Twisted forms}\label{twisted.forms.sect}

    Many properties are preserved by Tate-Shafarevich twists, even more by translation-twists.
Below is a list of known results, with sketches of proofs and some references.

    Starting with (\ref{vhs.cons.say}) we work over $\c$.
       
       \begin{say}[Properties preserved by twists I]\label{preserve.say}
         Assume that $\pi_i:X_i\to S$ are
           twisted forms of each other for $i=1,2$.

         (\ref{preserve.say}.1)  $X_1, X_2$  have the same singularities, up to \'etale isomorphism. In particular, if $X_1$ is smooth then so is $X_2$.
         If $X_1$ has normal, rational, canonical or terminal singularities, then the same holds for $X_2$.

         (\ref{preserve.say}.2) $\q$-factoriality is not an \'etale-local property.

         (\ref{preserve.say}.3)  If  $X_1$ has canonical singularities and $K_{X}$ is numerically $\pi_1$-trivial or $\pi_1$-nef, then the same hold for $X_2$.

         (\ref{preserve.say}.4)  If $X_i\to S$ are Abelian fiber spaces, then
         they have the same Picard number.
       \end{say}

       \begin{say}[Variations of Hodge structures]\label{vhs.cons.say}
          Let $\pi:X\to S$ be a projective morphism with  connected fibers, $X$ smooth and $S$ normal. Let $S^\circ$ be a dense open subset such that
          $\pi^\circ:X^\circ\to S^\circ$ is smooth. Thus we get integral
          variations of Hodge structures
          $$
          R^i\pi^\circ_*\z_{X^\circ}\to R^i\pi^\circ_*\c_{X^\circ}.
          \eqno{(\ref{vhs.cons.say}.1)}
          $$
          The following  properties are unchanged by  twisting.
          \begin{enumerate}\setcounter{enumi}{1}
          \item The $J$-invariant, mapping a point $s\in S^\circ$ to the Hodge structure  on $H^i(X_s, \z)$. 
          \item The local monodromy representations around the irreducible components of            $S\setminus S^\circ$.  
          \item  If $U\to S$ is \'etale, then the image of  $\pi_1(U)\to \pi_1(S^\circ)$ has finite index. 
            Thus if the $X_i\to S$ are twisted forms of each other, then the monodromy representations $\pi_1(S^\circ, s)\to \GL\bigl(H^i(F_s, \z)\bigr)$,
             agree on a finite index subgroup of
             $\pi_1(S^\circ, s)$.
            \end{enumerate}
\end{say}

       \begin{say}[Canonical bundle formula]\label{can.b.form.say}
          Assume that $X$ has canonical singularities and   $K_{X}$ is numerically $\pi$-trivial. The general Kodaira formula writes
              $K_{X}\simq \pi^*(K_S+J+B)$ where $J$ depends only on the
          J-invariant map (\ref{vhs.cons.say}.2), and $B$  depends only on the \'etale-local structure of $X\to S$  over  the codimension 1  points of  $S$; see \cite{k-adj} for an overview.

           Thus if the $X_i\to S$ are twisted forms of each other, then, 
          up to $\q$-linear equivalence, the  $K_{X_i}$ are pull-backs of the same $\q$-divisor class.

         As a special case, if $K_{X_1}\simq 0$ then
          $K_{X_2}\simq 0$.
\end{say}

       \begin{say}[Monodromy  preserving twists]\label{preserve.J.say}
         Assume that  $X$ has canonical singularities.
        By \cite[p.171]{k-hdi2}, the sheaves  $R^i\pi_* \omega_{X}$
         are determined by the variation of the integral Hodge structures  (\ref{vhs.cons.say}.1). The same holds for
        the sheaves  $R^i\pi_* \o_{X}$ by \cite[3.8]{k-hdi1}.

         Thus if the $X_i\to S$ are twisted forms of each other, and   the monodromy representations $\pi_1(S^\circ, s)\to \GL\bigl(H^i(F_s, \z)\bigr)$
             agree on the whole
             $\pi_1(S^\circ, s)$, then we have the following.

         (\ref{preserve.J.say}.1)   There are  natural isomorphisms
         $$
      R^i(\pi_1)_* \omega_{X_1}\cong R^i(\pi_2)_* \omega_{X_2},
          \qtq{and}
         R^i(\pi_1)_* \o_{X_1}\cong R^i(\pi_2)_* \o_{X_2}.
         $$

         (\ref{preserve.J.say}.2) Assume that the  $X_i$ are proper. Using using \cite[p.172]{k-hdi2}, 
         there are  natural isomorphisms
                $$
         H^i\bigl(X_1, \omega_{X_1}\bigr)\cong
         H^i\bigl(X_2, \omega_{X_2}\bigr)
         \qtq{and}
         H^i\bigl(X_1, \o_{X_1}\bigr)\cong
         H^i\bigl(X_2, \o_{X_2}\bigr).
         $$
         
         (\ref{preserve.J.say}.3) Assume that the  $X_i$ are proper. If  
$K_{X_1}\sim 0$, then $K_{X_2}\sim 0$. Indeed,  $K_{X_2}\simq 0$ by (\ref{preserve.say}.4) and it has a section by  (\ref{preserve.J.say}.2), so it is trivial. 
  \end{say}

       \begin{say}[Translation-twists]\label{preserve.translate.say} Let $\pi_i:X_i\to S$ be translation-twists of each other. 

         (\ref{preserve.translate.say}.1) There are natural isomorphisms
         $$
         R^i(\pi_1)_*\z_{X_1}\cong  R^i(\pi_2)_*\z_{X_2} \qtq{for every $i$.}
         $$
           In particular, translation-twists  are monodromy  preserving.

           Indeed, the two sheaves in (\ref{preserve.translate.say}.1) are locally the same. Using the notation of (\ref{et.tw.defn}.1), the patchings at a point $s\in S$ differ by the action of $\phi_{ij}|_s$  on   $H^i(X_s, \z_{X_s})$.
           Since  $\phi_{ij}|_s$ is in   the connected component of $\auts(X_s)$ by (\ref{tors.aut.say}.1),
           it acts trivially on the cohomology groups.

           Note however, that the differentials in the Leray spectral sequence computing $H^m(X_j, \z)$ can be different,
thus the fundamental group and the integral cohomology groups of the $X_i$  can be different;
see   (\ref{enr.surf.exmp}).

            This is similar to the situation encountered for minimal models.  Different minimal models $X^{\rm m}$ of a 3-fold $X$ have isomorphic  cohomology groups
           $H^i(X^{\rm m}, \z)$, but usually the cohomology rings  $H^*(X^{\rm m}, \z)$  are  not  isomorphic to each other; see \cite[3.2.2]{k-etc}.

            \smallskip
 (\ref{preserve.translate.say}.2) For smooth $X_i$, there are natural isomorphisms of the exterior algebras
         $$
         \oplus_i(\pi_1)_*\Omega^i_{X_1}\cong  \oplus_i(\pi_2)_*\Omega^i_{X_2}.
         $$
                    To see this, note that on a smooth fiber $X_z$ we have an exact sequence
            $$
            0\to \pi^*\Omega_Z\otimes k(z)\to \Omega_X|_{X_z}\to \Omega_{X_z}\to 0.
            $$
            Translation action is trivial on $\pi^*\Omega_Z$ and also on the global sections of $\Omega_{X_z}^{\otimes m}$ for any $m$. Thus a translation-isomorphism does not change
            the gluing data of the sheaves $ \oplus_i(\pi_j)_*\Omega^i_{X_j}$
            on a dense open set. Since the sheaves 
            $ \oplus_i(\pi_j)_*\Omega^i_{X_j}$ are torsion free, they are isomorphic to each other.

            If the $X_i$ are singular, similar results hold for the reflexive hulls $\Omega^{[i]}$. 

               \smallskip
               (\ref{preserve.translate.say}.3)    As a consequence, if $X_1$ is
               holomorphic
               symplectic, then so is $X_2$.

       \end{say}


       \newcommand{\etalchar}[1]{$^{#1}$}
\def\cprime{$'$} \def\cprime{$'$} \def\cprime{$'$} \def\cprime{$'$}
  \def\cprime{$'$} \def\dbar{\leavevmode\hbox to 0pt{\hskip.2ex
  \accent"16\hss}d} \def\cprime{$'$} \def\cprime{$'$}
  \def\polhk#1{\setbox0=\hbox{#1}{\ooalign{\hidewidth
  \lower1.5ex\hbox{`}\hidewidth\crcr\unhbox0}}} \def\cprime{$'$}
  \def\cprime{$'$} \def\cprime{$'$} \def\cprime{$'$}
  \def\polhk#1{\setbox0=\hbox{#1}{\ooalign{\hidewidth
  \lower1.5ex\hbox{`}\hidewidth\crcr\unhbox0}}} \def\cdprime{$''$}
  \def\cprime{$'$} \def\cprime{$'$} \def\cprime{$'$} \def\cprime{$'$}
\providecommand{\bysame}{\leavevmode\hbox to3em{\hrulefill}\thinspace}
\providecommand{\MR}{\relax\ifhmode\unskip\space\fi MR }
\providecommand{\MRhref}[2]{%
  \href{http://www.ams.org/mathscinet-getitem?mr=#1}{#2}
}
\providecommand{\href}[2]{#2}

 \bigskip

  Princeton University, Princeton NJ 08544-1000, \

  \email{kollar@math.princeton.edu}

\end{document}